\documentclass[10pt,a4paper,twoside,leqno]{article}
\usepackage{amsfonts,amssymb}
\usepackage{eufrak}
\usepackage{amscd}

\newtheorem{defn}{Definition}

\vspace{0.5cm}

 4
\font\ebf=cmbx8
\font\erm=cmr8

\usepackage{graphicx}

\begin{document}

\thispagestyle{empty}

\noindent {\bf Prefab posets` Whitney numbers }

\vspace{0.7cm} {\it A. Krzysztof Kwa\'sniewski}

\vspace{0.2cm}

{\erm High School of Mathematics and Applied Informatics}

{\erm  Kamienna 17, PL-15-021 Bia\l ystok, Poland}

{\erm e-mail: kwandr@wp.pl

\vspace{0.2cm}

\noindent {\ebf Summary}

\noindent {\small We introduce a natural partial order   $\leq$ in
structurally natural  finite subsets the cobweb prefabs sets
recently constructed by the present author. Whitney numbers of the
second kind of the corresponding subposet which constitute
Stirling-like numbers` triangular array - are then calculated and
the explicit formula for them  is provided. Next - in the second
construction - we endow the set sums of  prefabiants with such an
another partial order that their their Bell-like numbers include
Fibonacci triad sequences introduced recently by the present
author  in order to extend famous relation between binomial Newton
coefficients and Fibonacci numbers  onto the infinity of their
relatives among which there are also the Fibonacci triad sequences
and binomial-like coefficients (incidence coefficients included).

\vspace{0.2cm}

 AMS Classification Numbers: 05C20, 11C08, 17B56 .

\vspace{0.2cm}


\vspace{0.2cm}

\noindent presentation  (November $2005$) at the Gian-Carlo Rota
Polish Seminar\\
$http://ii.uwb.edu.pl/akk/index.html$

\vspace{0.2cm}

\noindent \textit{to appear in}  Bull. Soc. Sci. Lett. Lodz.
October 2005.

\vspace{0.3cm}

\section{Introduction} The clue algebraic concept of combinatorics - prefab
(with associative and commutative composition) was introduced in
[1], see also [2,3]. In [4] the present author hadconstructed a
new broader class of prefab`s notion extending combinatorial
structure based on the so called cobweb posets  (see Section 1.
[4] for the definition of a cobweb poset as well as a
combinatorial interpretation of its characteristic binomial-type
coefficients - for example- fibonomial ones [5,6]).\\
Here  we introduce two natural partial orders: one  $\leq$ in
grading-natural subsets of cobweb`s prefabs sets [4] and in the
second proposal we endow the set sums of prefabiants with such
another partial order that one  arrives at Bell-like numbers
including Fibonacci triad sequences introduced by the present
author in [7].

\vspace{1mm}

\section{Prefab based posets and their Whitney numbers.}
\vspace{1mm} Let the family $S$ of combinatorial objects
($prefabiants$) consists of all layers   $\langle\Phi_k
\rightarrow \Phi_n \rangle,\quad k<n,\quad k,n \in
N\cup\{0\}\equiv Z_\geq$ and an empty prefabiant $i$.

\noindent The set $\wp$ of prime objects consists of all
sub-posets $\langle\Phi_0 \rightarrow \Phi_m \rangle$  i.e. all
$P_m$`s $m \in N\cup\{0\}\equiv Z_\geq$ constitute from now on a
family of prime $prefabiants$ [4].

\noindent Layer is considered here to be the set of  all
max-disjoint isomorphic copies (iso-copies) of $P_{n-k}= P_{m}$
[4].

\noindent As a matter of illustration we quote after [4] examples
of cobweb posets` Hasse Diagrams so that the layers become
visualized.

\begin{center}
\includegraphics[width=75mm]{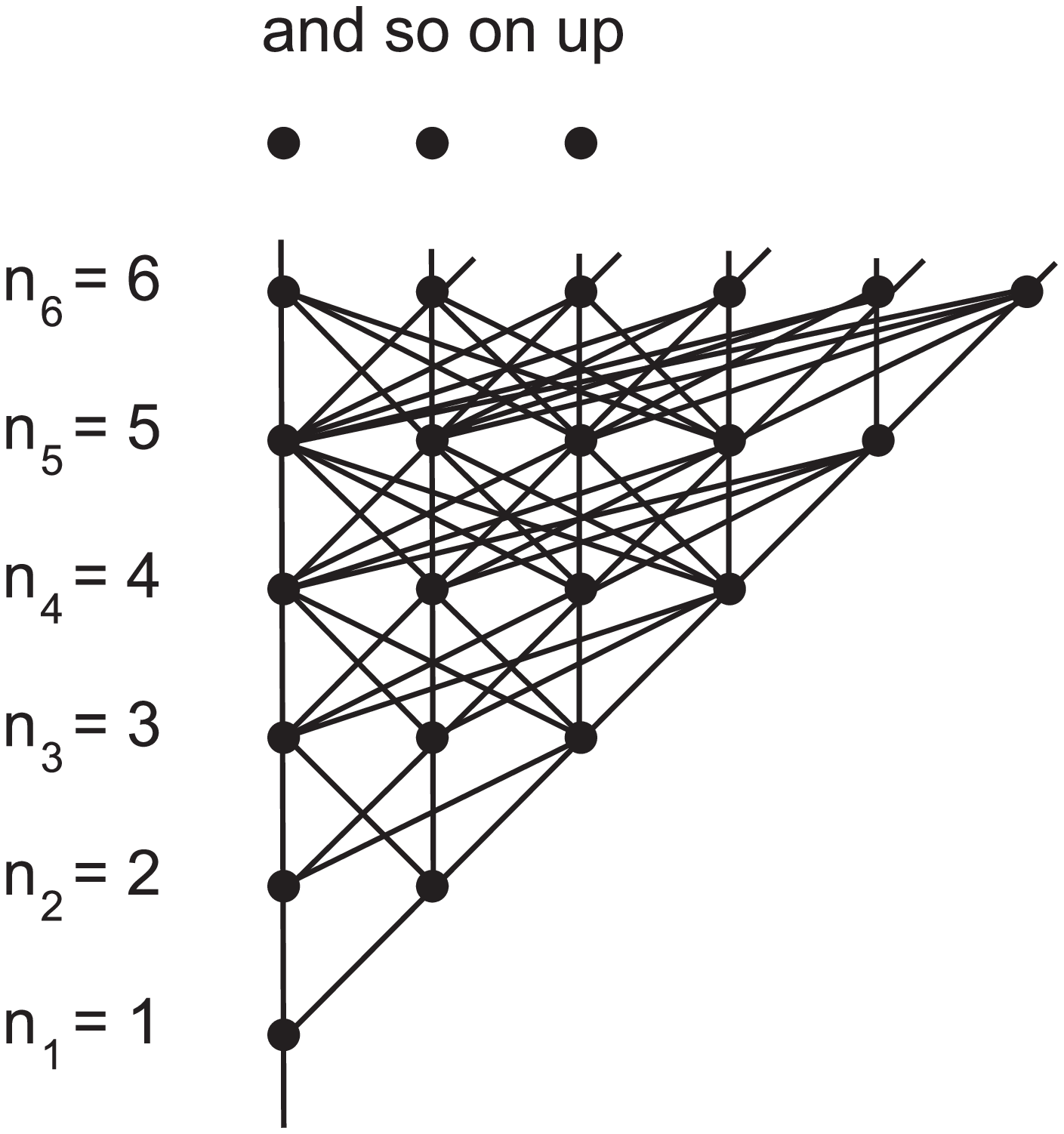}

\vspace{2mm}

\noindent {\small Fig.1. Display of Natural numbers` cobweb
poset.}
\end{center}

\vspace{2mm}

\begin{center}

\includegraphics[width=75mm]{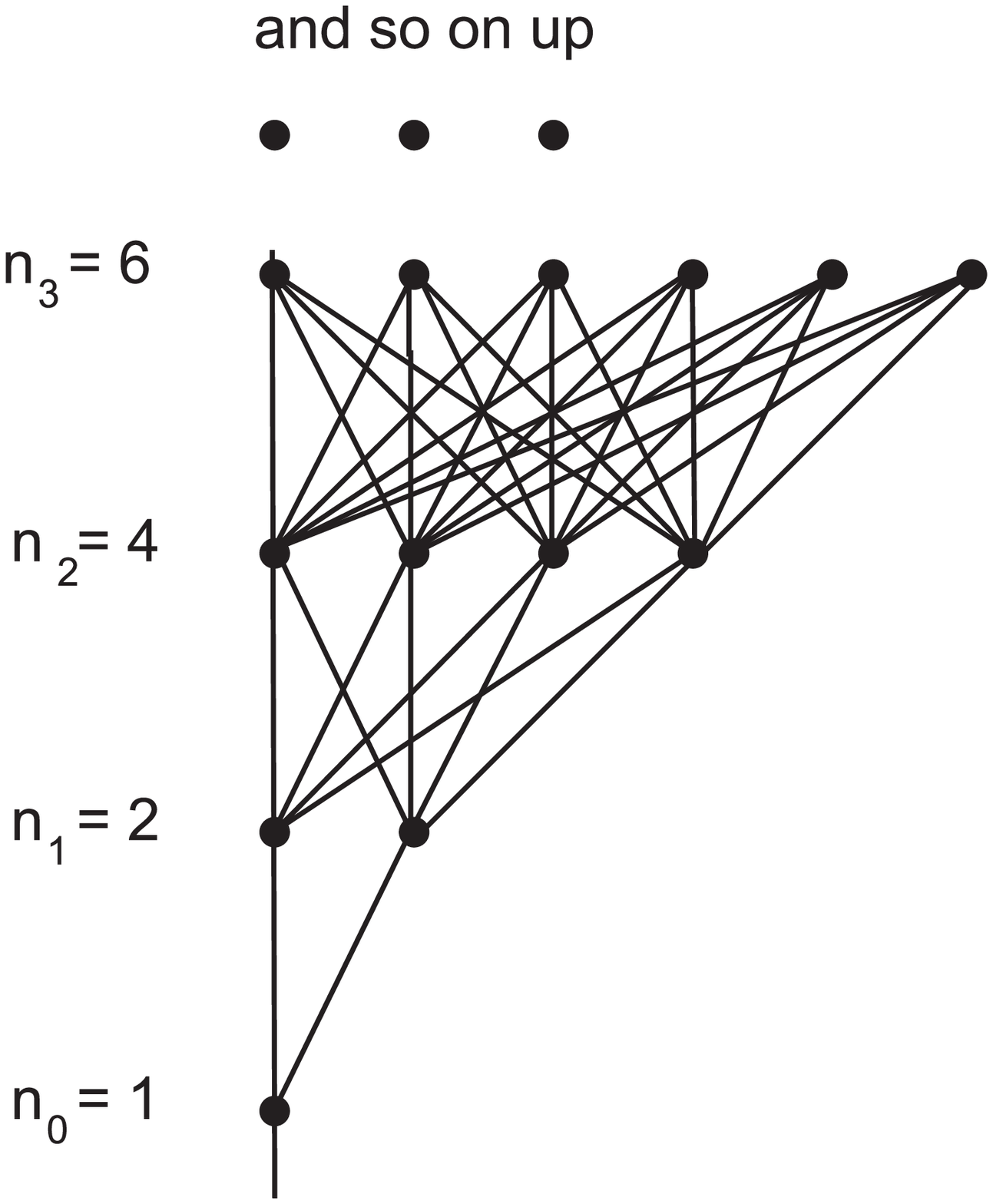}

\noindent {\small Fig.2. Display of Even Natural numbers  $\cup
\{1\}$ - cobweb poset.}

\end{center}

\vspace{2mm}

\begin{center}

\includegraphics[width=75mm]{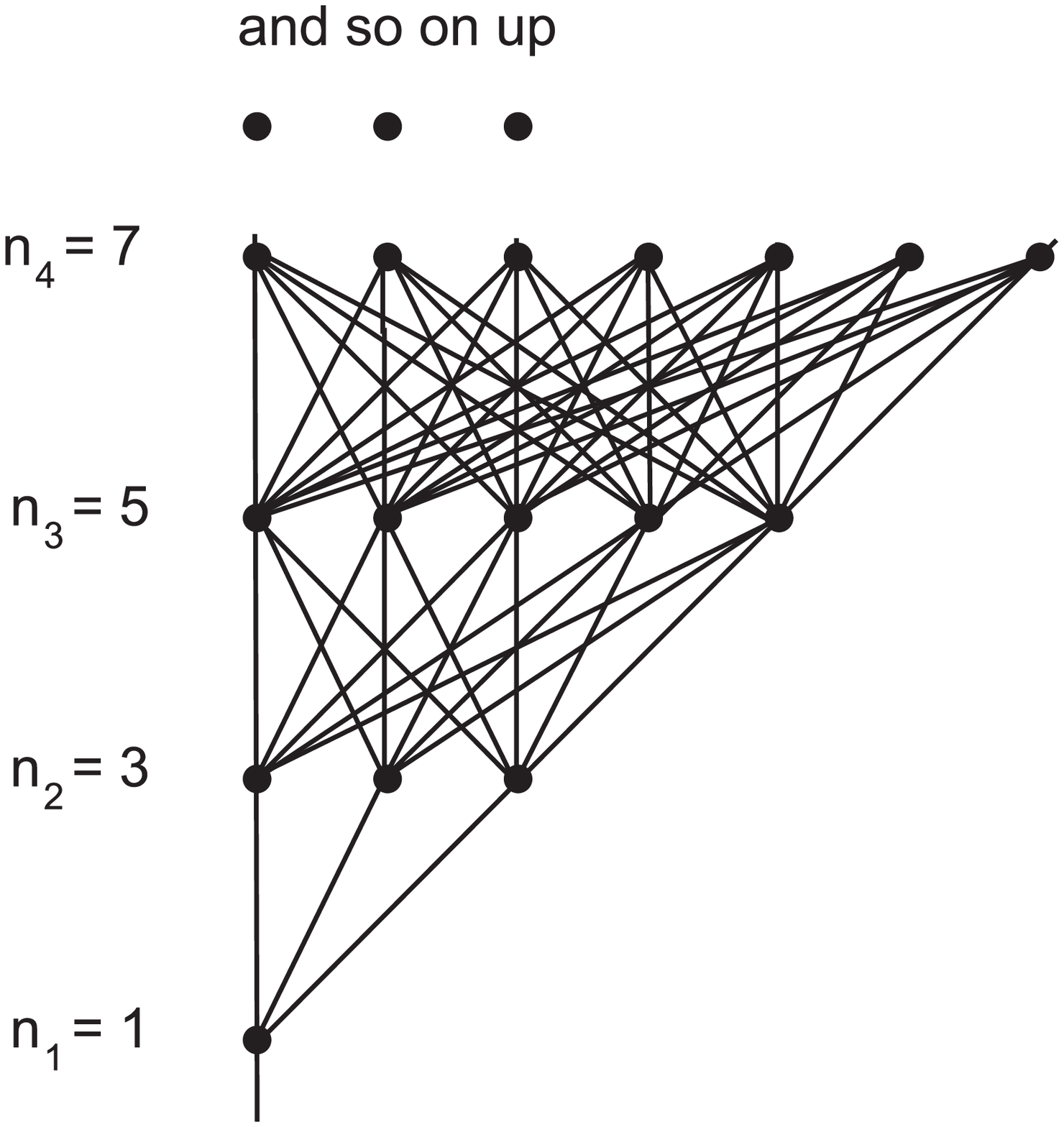}

\vspace{2mm}

\noindent {\small Fig3. Display of Odd natural numbers` cobweb
poset.}
\end{center}

\vspace{2mm}

\begin{center}

\includegraphics[width=75mm]{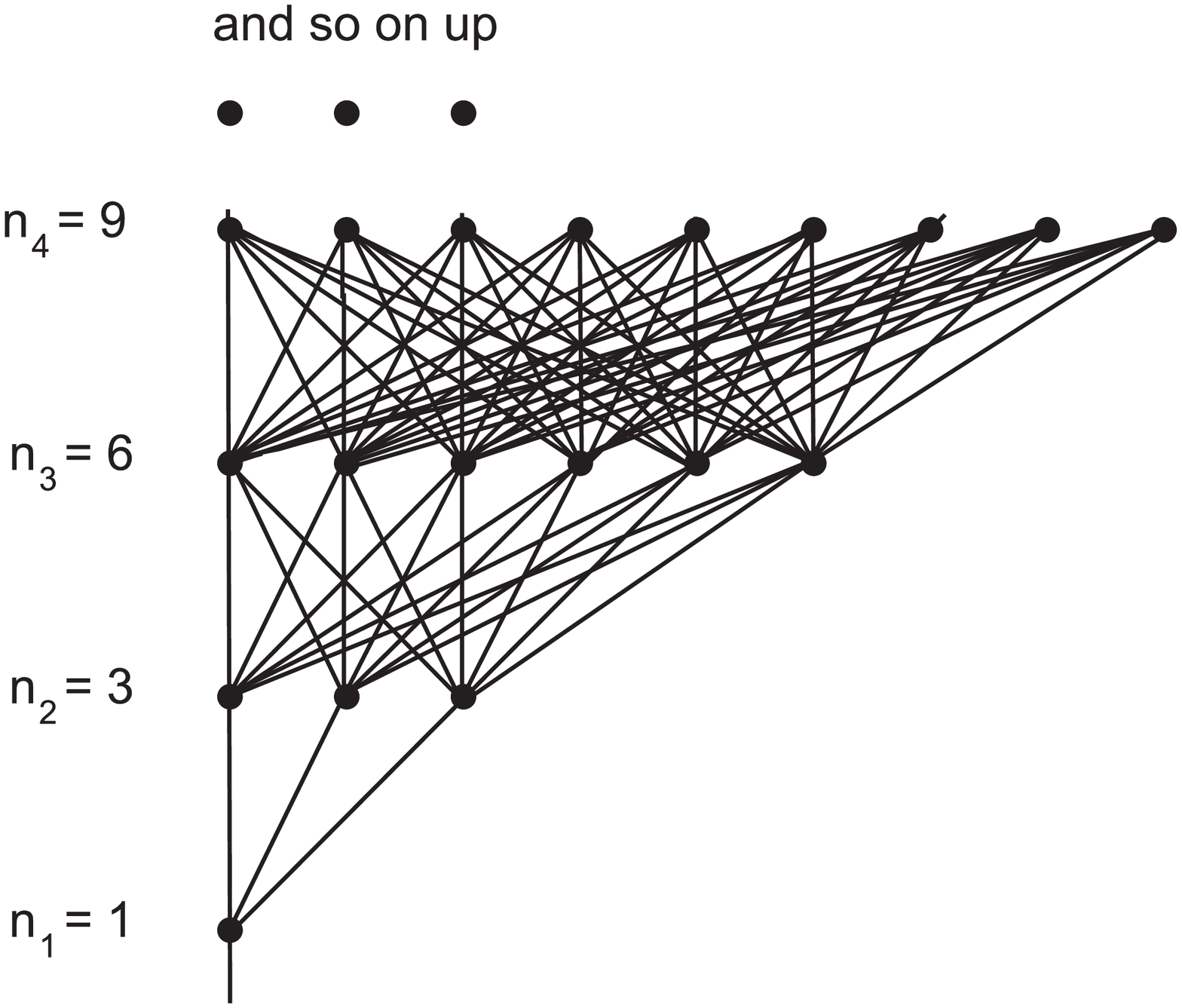}

\vspace{2mm}

\noindent {\small Fig.4. Display of divisible by 3 natural numbers
$\cup \{1\}$ - cobweb poset.}

\end{center}

\vspace{2mm}

\begin{center}

\includegraphics[width=75mm]{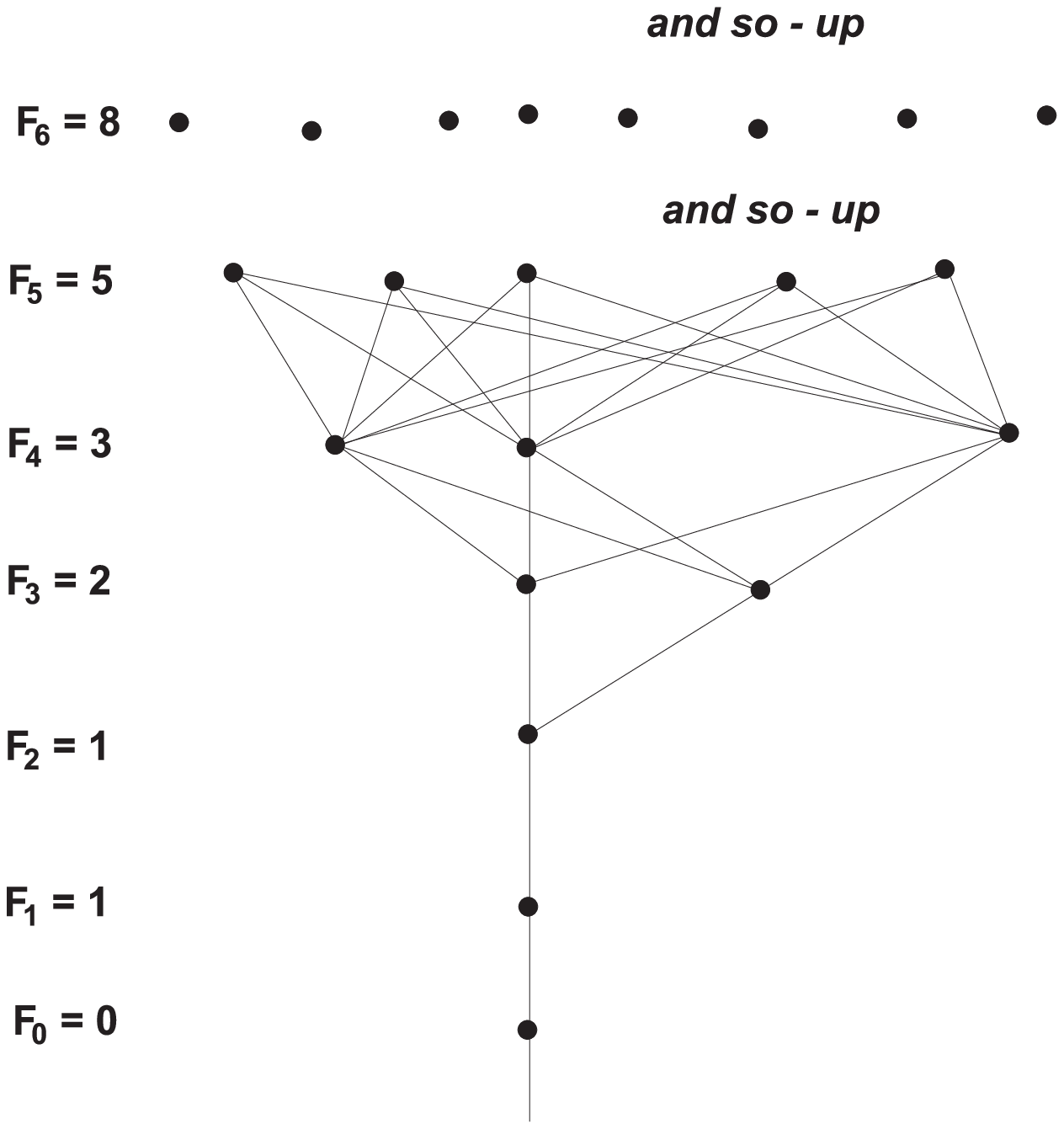}

\vspace{2mm}

\noindent {\small Fig.5. Display of Fibonacci numbers` cobweb
poset.}
\end{center}

\vspace{0.1cm}

\noindent Consider then now  the partially ordered family $S$ of
these layers considered to be sets of  all max-disjoint isomorphic
copies (iso-copies) of  prime prefabiants $P_{m}= P_{n-k}$
displayed by Fig 1. - Fig.5. above [4]. Let us define  in $S$ the
partial order relation as follows.

\begin{defn}
$$
\langle \Phi_k \rightarrow \Phi_n \rangle \leq
\langle\Phi_{k^\clubsuit} \rightarrow \Phi_{n^\clubsuit}\rangle
\quad \equiv\quad k \leq k{^\clubsuit} \quad\wedge\quad n\leq
n^{\clubsuit}.
$$
\end{defn}
For convenience reasons we shall also adopt and use the following
notation: $$ \langle\Phi_k \rightarrow \Phi_n \rangle = p_{k,n}.
$$
The interval $[p_{k,n}, p_{{k^\clubsuit},{n^\clubsuit}}]$  is of
course a subposet of $\langle\ S , \leq \rangle$. We shall
consider in what follows the subposet  $\langle\ P_{k,n} , \leq
 \rangle$ where
$$ P_{k,n}= [p_{o,o} ,p_{k,n}].$$

\noindent {\bf Observation 1.} The size $|P_{k,n}|$ of $P_{k,n}$ $
=|\{\langle l,m \rangle ,\quad 0 \leq l \leq k \quad\wedge \quad 0
\leq m \leq n \quad\wedge \quad k\leq n \}| = (n-k)(k + 1) +
\frac{k(k+1)}{2}.$

\vspace{2mm}

\noindent Proof: Obvious. Just draw the picture $\{\langle l,m
\rangle ,\quad 0 \leq l \leq k \quad\wedge \quad 0 \leq m \leq n
\quad\wedge \quad k\leq n \}$ of $P_{k,n}$` grid.

\vspace{1mm}

\noindent {\bf Observation 2.} The number of maximal chains in
$\langle\ P_{k,n} , \leq \rangle$ is equal to the number $d(k,n)$
of $0$ - dominated strings of binary i.e.  $0's$  and $1's$
sequences

$$d(k,n) = \frac{n+1-k}{n} \left( \begin{array}{c}
{k+n}\\n\end{array}\right). $$

\vspace{2mm}

\noindent Proof. The number we are looking for equals to the
number of minimal walk-paths in Manhattan grid [8] $[k\times n]$
restricted by the condition $ k \leq n $  i.e. it equals to the
number of $0$ - dominated strings of  $0's$  and  $1's$ sequences.

\vspace{2mm}

\noindent Recall that $( d(k,n) )$  infinite matrix`s  diagonal
elements are equal to the \textbf{Catalan} numbers $C(n)$

$$C(n) = \frac{1}{n} \left( \begin{array}{c}
{2n}\\n\end{array}\right). $$ as the Catalan numbers count the
number of  $ 0$ - dominated strings of   $0's$  and $1's$   with
equal number of $0's$  and $1's$ . Recall that  a $0$ - dominated
string of length  $n$  is such a string  that the first  $k$
digits of the string contain at least as many  $0's$ as $1's$ for
$k = 1, . . . , n $  i.e.  $0`s$ prevail  in appearance, dominate
$1`s$  from the left to the right end of the string.  $ 0$ -
dominated strings correspond bijectively  to minimal bottom - left
corner to the right upper corner paths in an integer grid $Z_\geq
\times Z_\geq$ rectangle part called Manhattan [8] with the
restriction imposed on those minimal paths to obey the "safety"
condition $k \leq n$ .

\vspace{2mm}

\noindent \textbf{Comment 1.} Observation 2. equips the poset
$\langle\ P_{k,n} , \leq \rangle$ with clear cut combinatorial
meaning.

\vspace{2mm}

\noindent The poset $\langle\ P_{k,n} , \leq \rangle$ is naturally
graded. $\langle\ P_{k,n} , \leq \rangle$ poset`s maximal chains
are of all of equal size (Dedekind property) therefore the rang
function is defined.

\vspace{2mm}

\noindent {\bf Observation 3.} The rang $r(P_{k,n})$  of $P_{k,n}$
= number of elements in maximal chains $P_{k,n}$ minus one
$=k+n-1$. \noindent The rang $r(p_{l,m})$ is defined accordingly:
$r(p_{l,m})= l+m-1.$

\vspace{2mm}

\noindent Proof: obvious. Just draw the picture $\{\langle l,m
\rangle ,\quad 0 \leq l \leq k \quad\wedge \quad 0 \leq m \leq n
\quad\wedge \quad k\leq n \}$  of $P_{k,n}$` grid and note that
maximal means paths without at a slant edges.

\vspace{2mm}

\noindent Accordingly Whitney numbers $W_k(P_{l,m})$  of the
second kind are defined as follows (association: $ n
\leftrightarrow  \langle l,m \rangle $)

\vspace{2mm}

\begin{defn}
$$
W_k(P_{l,m})= \sum_{\pi\in P_{l,m}, r(\pi)=k} 1 \quad\equiv\quad
S(k,\langle l,m \rangle).
$$
\end{defn}

\vspace{2mm}

\noindent Here now and afterwards we identify  $W_k(P_{l,m})$ with
$S(k,\langle l,m \rangle)$ called and  viewed at  as Stirling -
like numbers of the second kind  of the naturally graded poset
$\langle\ P_{k,n} , \leq \rangle$ - note the association: $ n
\leftrightarrow  \langle l,m \rangle $.

\vspace{2mm}

\noindent \textbf{Right now challenge problems.  I.}

\vspace{2mm}

\noindent \textbf{I.} Let us define now Whitney numbers
$w_k(P_{l,m})$  of the first kind as follows (association: $ n
\leftrightarrow \langle l,m \rangle $.  Note the text-book
notation for M\"{o}bius function $\mu$)

\vspace{2mm}

\begin{defn}
$$
w_k(P_{l,m})= \sum_{\pi\in P_{l,m}, r(\pi)=k} \mu(0,\pi) \equiv
s(k,\langle l,m \rangle).
$$
\end{defn}

\vspace{0,2cm}

\vspace{2mm}

\noindent Here now and afterwards we identify  $w_k(P_{l,m})$ with
$s(k,\langle l,m \rangle)$  called and viewed at as Stirling -
like numbers of the first kind  of the poset $\langle\ P_{k,n} ,
\leq \rangle$ - note the association: $ n \leftrightarrow \langle
l,m \rangle $.

\noindent \textbf{Problem 1} Find an explicit expression for
$$
w_k(P_{l,m})\equiv s(k,\langle l,m \rangle) = ?
$$
and
$$
W_k(P_{l,m})\equiv S(k,\langle l,m \rangle) = ?
$$
Occasionally note that $S(k,\langle l,m \rangle)$ equals to the
number of the grid points counted at a slant (from the up-left to
the right-down) accordingly to the $l+m = k$ requirement.

\vspace{2mm}

\noindent \textbf{Problem 2} Find the recurrence relations for
$$
w_k(P_{l,m})\equiv s(k,\langle l,m \rangle \quad and\quad
W_k(P_{l,m})\equiv S(k,\langle l,m \rangle).
$$
We define now (note the association: $ n \leftrightarrow \langle
l,m \rangle $)  the corresponding Bell-like numbers
$$ B(\langle l,m \rangle)$$
of the naturally graded poset $\langle\ P_{k,n} , \leq \rangle$ as
follows.

\begin{defn}
$$
B(\langle l,m \rangle)= \sum_{k}^{l+m}S(k,\langle l,m \rangle).
$$
\end{defn}

\vspace{2mm}

\noindent {\bf Observation 4.}
$$ B(\langle l,m \rangle)=
|P_{l,m}| = \frac{k(k+1)}{2} + (n-k)(k+1).
$$
Proof: Just draw the picture $\{\langle l,m \rangle ,\quad 0 \leq
l \leq k \quad\wedge \quad 0 \leq m \leq n \quad\wedge \quad k\leq
n \}$ of $P_{k,n}$` grid and note that $S(k,\langle l,m \rangle)$
equals to the number of the grid points counted at a slant (from
the up-left to the right-down) accordingly to the $l+m = k$
requirement. Summing them up over all gives the size of $P_{k,n}$.

\vspace{2mm} \noindent \textbf{Comment 2.} Observation 4. equips
the poset`s $\langle\ P_{k,n} , \leq \rangle$ Bell-like numbers
$B(\langle l,m \rangle)$ with clear cut combinatorial meaning.

\section{Set Sums of prefabiants` posets and their Whitney numbers.}
\vspace{1mm} In this part we consider prefabiants` set sums with
an appropriate another partial order so as to arrive at Bell-like
numbers including Fibonacci triad sequences introduced recently by
the present author in [7] - see also [9].

\noindent Let $F$ be any   \textit{ "GCD-morphic"} sequence [4].
This means that $GCD[F_n,F_m] = F_{GCD[n,m]}$ where $GCD$ stays
for Greatest Common Divisor operator. We define the finite partial
ordered set $P(n,F)$ as the set of \textbf{prime} prefabiants
$P_l$ given by the sum below.

\begin{defn}
$$ P(n,F) =\bigcup_{0\leq p}\langle \Phi_p \rightarrow \Phi_{n-p} \rangle = \bigcup_{0\leq l}P_{n-l}$$
\end{defn}
with the partial order relation defined  for $n-2l\leq 0$
according to
\begin{defn}
$$
P_l \leq P_{\hat l} \quad \equiv \quad l\leq \hat l,\quad P_{\hat
l}, P_l\in \langle \Phi_l \rightarrow \Phi_{n-l} \rangle.
$$
\end{defn}
Recall that \textbf{rang of} $P_l$ \textbf{is} $l$. Note that
$\langle \Phi_l \rightarrow \Phi_{n-l} \rangle = \emptyset$ for
$n-2l \leq 0$. The Whitney numbers of the second kind are
introduce accordingly.

\vspace{2mm}

\begin{defn}
$$
W_k(P_{n,F})= \sum_{\pi\in P_{n,F}, r(\pi)=k} \equiv S(n,k,F).
$$
\end{defn}
Right from the definitions above we infer that: (recall that
\textbf{rang of} $P_l$ \textbf{is} $l$.)

\vspace{2mm}

\noindent \textbf{Observation 5.}
$$
W_k(P_{n,F})= \sum_{\pi\in P_{n,F}, r(\pi)=k} \equiv S(k,n-k,F)=
\left( \begin{array}{c} {n-k}\\k\end{array}\right)_{F}.
$$
Here now and afterwards we identify $W_k(P_{n,F})=S(n,k,F)$ viewed
at and called as Stirling - like numbers of the second kind of the
$P$ defined in [4]. $P$ by construction (see Figures above)
displays self-similarity property with respect to its prime
prefabiants sub- posets $P_n = P(n,F)$.

\vspace{2mm}

\noindent \textbf{Right now challenge problems. II.}

\noindent We repeat with obvious replacements of corresponding
symbols, names and definitions the same problems as in "Right now
challenge problems. I".

\noindent Here now consequently - for any $GCD$-morphic  sequence
$F$  (see: [4])  we define the corresponding Bell-like numbers $
B_n(F)$ of the poset $P(n,F)$ as follows.
\begin{defn}
$$B_n(F)= \sum_{k\geq0}S(n,k,F).$$
\end{defn}
Due to the investigation in [9,7] we have right now at our
disposal all corresponding results of [7,9] as the following
identification with  special case of $\langle
\alpha,\beta,\gamma\rangle$ - Fibonacci sequence $\langle
F_n^{[\alpha,\beta,\gamma]}\rangle_{n\geq 0}$ defined in [7]
holds.

\vspace{2mm}

\noindent {\bf Observation 6.}
$$B_n(F)\equiv F_{n+1}^{[\alpha=0,\beta=0,\gamma=0].}$$
Proof: See the Definition 2.2. from [7]. Compare also with the
special case of  formula (6) in [9].

\vspace{2mm}

\noindent\textbf{ Recurrence relations.} Recurrence relations for
$\langle \alpha,\beta,\gamma\rangle$ - Fibonacci sequences
$F_n^{[\alpha,\beta=,\gamma]}$ are to be found in [7] - formula
(9). Compare also with the special case formula (7) in [9].

\vspace{2mm}

\noindent\ \textbf{Closing-Opening Remark.} The study of further
properties of these Bell-like numbers as well as the study of
consequences of these identifications for the domain of the
widespread data types [7] and perhaps for eventual new dynamical
data types we leave for the possibly coming future. Examples of
special cases - a bunch of them - one finds in [7] containing [9]
as a special case. As seen from the identification Observation 6.
the special cases of $\langle \alpha,\beta,\gamma\rangle$ -
Fibonacci sequences $F_n^{[\alpha,\beta,\gamma]}$  gain
\textbf{additional} with respect to [9,7] combinatorial
interpretation in terms Bell-like numbers as sums over  rang $=k$
parts of the poset i.e. just sums of Whitney numbers of the poset
$P(n,F)$. This adjective \textit{"additional"} shines brightly
over Newton binomial connection constants between bases
$\langle(x-1)^k\rangle_{k\geq0}$ and $\langle x^n\rangle_{n\geq0}$
as these are Whitney numbers of the numbers from $[n]$ chain i.e.
Whitney numbers of the poset $\langle [n], \leq \rangle.$ For
other elementary "shining brightly"  examples see Joni , Rota and
Sagan excellent presentation in [10].

\vspace{2mm}

\noindent  \textbf{Acknowledgements}

\noindent Discussions with Participants of Gian-Carlo Rota Polish
Seminar  on all related topics \\
$http://ii.uwb.edu.pl/akk/index.html$  - are appreciated with
pleasure.

\begin
{thebibliography}{99}
\parskip 0pt

\bibitem{1}
E. Bender, J. Goldman   {\it Enumerative uses of generating
functions} , Indiana Univ. Math.J. {\bf 20} 1971), 753-765.

\bibitem{2}
D. Foata and M. Sch"utzenberger, Th'eorie g'eometrique des
polynomes euleriens, (Lecture Notes in Math., No. 138).
Springer-Verlag, Berlin and New York, 1970.

\bibitem{3}
A. Nijenhuis and H. S. Wilf, Combinatorial Algorithms, 2nd ed.,
Academic Press, New York, 1978.

\bibitem{5}
A. K.  Kwa\'sniewski, {\it Cobweb posets as noncommutative
prefabs} submitted for publication ArXiv : math.CO/0503286 (2005)

\bibitem{6}
A. K.  Kwa\'sniewski, {\it Information on combinatorial
interpretation of Fibonomial coefficients }   Bull. Soc. Sci.
Lett. Lodz Ser. Rech. Deform. 53, Ser. Rech.Deform. {\bf 42}
(2003), 39-41. ArXiv: math.CO/0402291   v1 18 Feb 2004

\bibitem{6}
A. K. Kwa\'sniewski, {\it The logarithmic Fib-binomial formula}
Advanced Stud. Contemp. Math. {\bf 9} No 1 (2004), 19-26. ArXiv:
math.CO/0406258 13 June 2004.

\bibitem{7}
A. K.  Kwa\'sniewski, {\it Fibonacci-triad sequences} Advan. Stud.
Contemp. Math. {\bf 9} (2) (2004),109-118.

\bibitem{8}
Z. Palka , A. Ruciñski {\it Lectures on Combinatorics}.I. WNT
Warsaw 1998 (\textit{in polish})

\bibitem{9}
A. K.  Kwa\'sniewski, {\it Fibonacci q-Gauss sequences} Advanced
Studies in Contemporary Mathematics {\bf 8} No 2 (2004), 121-124.
ArXive:  math.CO/0405591  31 May 2004.

\bibitem{10}
S.A. Joni ,G. C. Rota, B. Sagan {\it From sets to functions: three
elementary examples} Discrete Mathematics {\bf 37} (1981),
193-2002.

\end{thebibliography}



\end{document}